\documentclass[a4paper,10pt]{article}
\usepackage{amsfonts,amssymb,amsmath,amsthm}
\usepackage{mathrsfs}
\usepackage{latexsym}
\usepackage[hmargin=2cm,vmargin=2.5cm]{geometry}

\title{Eigenvalue gaps for the Laplacian on hypersurfaces of the
sphere}
\author{Demetrios A. Pliakis}

\begin{document}

\newcommand{\hmi}{\ensuremath{\overline{{\bf R}_+}}}
\newcommand{\mom}[1]{\ensuremath{{#1}:{\cal H}\rightarrow {\bf R}}}
\newcommand{\msm}[1]{\ensuremath{{#1}:{\cal S}\rightarrow {\bf R}}}
\newcommand{\sphe}[1]{\ensuremath{({\bf S}^{#1},e)}}
\newcommand{\kseks}{\begin{eqnarray}}
\newcommand{\tleks}{\end{eqnarray}}
\newcommand{\kseksw}{\begin{eqnarray*}}
\newcommand{\tleksw}{\end{eqnarray*}}
\newcommand{\lap}[1]{\ensuremath{\Delta_{#1}}}
\newcommand{\laps}{\ensuremath{\Delta_{{\cal S}}}}
\newcommand{\lapv}{\ensuremath{\Delta_{{\cal H}}}}
\newcommand{\ei}{\ensuremath{\lambda}}
\newcommand{\eps}{\ensuremath{\epsilon}}
\newcommand{\inm}[1]{\ensuremath{\int_{#1} }}
\newcommand{\Ha}[1]{\ensuremath{ {\cal H}^{#1} }}
\newcommand{\eaf}[1]{\ensuremath{\sup_{#1}}}
\newcommand{\mkf}[1]{\ensuremath{\inf_{#1}}}
\newcommand{\spe}[5]{\ensuremath{ \mbox{V}_{#1,#2}(#3,#4,#5)}}
\newcommand{\sta}[6]{\ensuremath{\mbox{C}_k(#1,#2,#3,#4,#5,#6)}}
\newcommand{\no}[2]{\ensuremath{||#1||_{#2}  }}
\newcommand{\cinf}[2]{\ensuremath{C_0^\infty({\bf R}^{#1}{#2})}}
\newcommand{\bge}{\ensuremath{B_\gamma(\epsilon_\gamma)}}
\newcommand{\ric}{\ensuremath{\mbox{Ric}}}
\newcommand{\R}{\ensuremath{\mbox{R}}}
\newcommand{\rie}{\ensuremath{\mbox{Rm}}}
\newcommand{\se}{\ensuremath{\mbox{K}}}
\newcommand{\xw}[1]{\ensuremath{{\bf R}^{#1}}}
\newcommand{\sfe}[1]{\ensuremath{{\bf S }^{#1}}}
\newcommand{\Sy}[1]{\ensuremath{\mbox{S}_{#1}}}
\newcommand{\hare}[1]{\ensuremath{|\iP|^{#1}}}
\newcommand{\kl}{\ensuremath{\nabla}}
\newcommand{\iq}{\ensuremath{\mbox{tr}}}
\newcommand{\es}{\ensuremath{\mbox{Hess}}}
\newcommand{\rs}{\ensuremath{\mbox{R}}}
\newcommand{\hard}[1]{\ensuremath{\left|\frac{\nabla{#1}}{{#1}}\right|^2}}
\newcommand{\fo}{\ensuremath{\mbox{supp}}}
\newcommand{\iP}{\ensuremath{\mbox{P}}}
\newcommand{\olwv}[1]{\ensuremath{\int_{{\cal H}_{#1}}}}
\newcommand{\dian}[1]{\ensuremath{\underline{#1}}}
\newcommand{\olws}{\ensuremath{\int_{{\cal S}} }}
\newtheorem{hc}{Theorem}[section]
\newtheorem{alf}{Theorem}[section]

\maketitle

\begin{abstract}

We provide lower estimates for the eigenvalues of the laplacian for 
hypersurfaces of the round sphere. 

\end{abstract}

\paragraph{Introduction} The laplacian acting on functions on a compact 
riemannian 
manifold exhibits a discrete spectrum of positive eigenvalues counted with
multiplicity:
\[
 0< \lambda_1\leq \lambda_2 \leq \dots \leq \lambda_n \leq \dots 
\]
We provide an estimate for the gaps between eigenvalues:
\[
 \lambda_{i+1}-\lambda_i \geq C_i>1 
\]
in the case of a convex analytic hypersurface 
\({\cal H}\subset \sphe{n+1},n>1\) 
in the round sphere with constant depending on its second fundamental form. 
The problem appears interesting in the construction of minimal hypersurfaces
with singularities and \({\cal H}\) is then a minimal hypersurface 
in \(\sphe{n+1}\).

\paragraph{The harmonic extension method}

We will apply the method of harmonic extension from \cite{cw} or \cite{sy}  
Let $f:{\cal H}\rightarrow {\bf R}$ be a function 
defined on a smooth hypersurface \({\cal H}\subset \sphe{n+1}\), dividing 
\(\sphe{n+1}\) 
in two  regions \({\cal S}_1,{\cal S}_2\) that we denote simply by
\({\cal S}\).  
Let  \(\msm{U}\) be the solution of  the following boundary value 
problem:
\[
 \laps U= 0,\,\,\,\, U|_{{\cal H}}=f
\]
Then we have that if $\dian{N}$ is chosen as the unit normal to \({\cal H}\) 
and  \(H\) its mean curvature then
\[
0=\laps U= \dian{N}\left(\dian{N}U\right)+\lapv f +nH\dian{N}(U)
\]
Weitzenb\"ock formula suggests that for a harmonic function
\[
 \laps|\kl U|^2 = 2 |\es U|^2+2n|\kl U|^2 \geq  2n|\kl U|^2 
\]
Hence we have that by an integration by parts 
\begin{equation}\label{BI}
 \olws \laps|\kl U|^2 =-\olwv{} w\left(4\lapv f+2nHw\right)-
2\olwv{} h(\kl f, \kl f) 
\end{equation}
where \(w(x)=\dian{N}(U)\).
This identity will be the basis of our considerations.  We assume that 
\(f\) posseses the unique conitnuation property on \({\cal H}\)
while $u$ is analytic in ${\cal S}$. Then $\dian{N}$ the normal vector 
field  to ${\cal H}$ that satisfies the unique continuation property 
(e.g. for a minimal variety) then the set the analytic weight function 
\(w(x)\)and 
\[
{\cal N}=\{ x\in {\cal H}/ w(x)=0\} 
\] 
is a recitifiable set of finite mulitplicity in the hypersurface
${\cal H}\subset \sphe{n+1}$. Actually Sard's lemma asserts that
for almost all \(\eta\in (0,\eta_0)\) the sets 
\[
{\cal N}_\eta=\{ x\in {\cal H}/ |w(x)|\leq \eta \} 
\] 
are smooth.  Then we have the subsets  
\[
{\cal H}_\pm=\{ x\in {\cal H}/ \pm w(x) >0\} 
\]
  where we need to choose $f$ in order to 
estimate the integrals $\olwv{\pm} w \lap{{\cal H}_\pm}f $.  
We will glue the two functions $f_\pm$ with a partition of unity $\chi_\pm$ 
and transition area a tubular neighbourhood ${\cal N}_\eta$ of ${\cal N}$ 
of thickness $\eta$:
after the \L ojasiewicz inequality for analytic functions, 
adapted to the case
of  an analytic submanifold of the sphere \(\sphe{n+1}\).
 Indeed we have  that the localizing functions satisfy 
\[
\mbox{supp}(\chi_\pm)\subset {\cal H}_\pm, \,\,\,\, \chi_\pm=1\,\,\,
 \mbox{in} \,\,\,\,  {\cal H}_\pm\setminus
{\cal N}_\eta={\cal H}_{\pm,\eta}, \,\,\,\, 
|\kl\chi_\pm|+\tau|\kl^2\chi_\pm|\leq\frac{C}{\tau} 
\] 
for the parameter \(0<\tau<1\) controlling the transition regions.
\paragraph{ The choice of $f$.} We select as
\[
f=\frac{(u^2+\alpha^2)}{(v^2+\beta^2)} 
\]
where $u,v$ is the eigenfunction corresponding to the eigenvalue
$\lambda,\mu$ and the $\pm$ sign depends on the ${\cal H}_\pm$ region. 
This function obeys the  differential equation on ${\cal H},
\delta=\lambda-\mu$
\[
\lapv f = 2\left[-\delta+\Pi\right]f
\]
where 
\[
\Pi=\frac{\lambda\alpha^2}{u+\alpha^2}-
\frac{\mu\beta^2}{u^2+\beta^2}-
\frac{4uv\kl u\cdot \kl v}{(u^2+\alpha^2)(v^2+\beta^2)}+
\frac{|\kl u|^2}{u^2+\alpha^2}
-\frac{\zeta|\kl v|^2}{v^2+\beta^2},\qquad \zeta =\frac{b^2-5v^2}{v^2+b^2}
,\qquad \zeta_0\leq \zeta \leq 1
\]
where \(\zeta_0\) is determined by \(\beta\).
\paragraph{\bf A. The ${\cal H}_{+,\eta}$ region.} In order to bound 
\(\Pi\) we appeal to Harnack and Berstein inequalities proved in \cite{pno}. 
Specifically, exhausting the region \({\cal H}_\eta\) through \(0<\theta<1\):
\[
{\cal H}_\eta=\bigcup_{j=1}^N {\cal G}_j,\qquad {\cal G}_j=
\{x\in {\cal H}/ \eta_0\theta^j\leq w(x)\leq 
\eta_0\theta^{j-1}\}
\]  
We have that for numerical constants  \(c_1,c_2>0,d\geq 1\) and 
near \(|u|\geq \varrho \) and \(|v|\geq \varrho\):
\[
\eaf{{\cal H}_\eta}|\kl u| \leq c_1 (\eta_0\theta^j)^{\frac{d}{n-2}}
\lambda^{dn+1}\varrho
\]
\[
\eaf{{\cal H}_\eta}|\kl v| \leq c_2 (\eta_0\theta^j)^{\frac{d}{n-2}}
\mu^{dn+1}\varrho
\]
If we assume that the hypersurface is real analytic - it holds for the
minimal case- then the \(\varrho\) 
tubular neighbourhood of the nodal sets of \(u,v\) varies as
\(\varrho\lambda^{\frac{n-1}{2}},\varrho\mu^{\frac{n-1}{2}}\).
and hence selecting \(\rho\) analogously we can make this arbitrarily small.
Young's inequality allows us to write that:
\[
 \Pi \geq 
\frac{\lambda\alpha^2+\left(1-2\epsilon\right)|\kl u|^2}{u^2+\alpha^2}
-\frac{\mu\beta^2+\left(\frac{2}{\epsilon}+\zeta\right)|\kl v|^2}{v^2+\beta^2}
\]
Carefull consideration of the preceding bounds allows us to 
select \(\alpha,\beta,\theta\) so that for \(n>2\)
\[
\theta=\frac{1}{(2\mu)^{\frac{n}{n-2}}\mu^{n(n-2)}},\quad 
\alpha=  2\sqrt{c_1} (\eta_0\theta^j)^{\frac{d}{n-2}} \lambda^{dn}\varrho,
\quad 
\beta= c_3\varrho
\]
The \(n=2\) case requires that 
\[
\theta=\frac{1}{\mu},\quad 
\alpha=  2\sqrt{c_1} (\eta_0\theta^j)^d \lambda^{2d}\varrho,
\quad 
\beta= c_3\varrho
\]
We conclude that 
\[
\Pi\geq \frac12\delta
\]
\paragraph{\bf B. The region ${\cal H}_{-,\eta}$}. In this part of the
bounding hypersurface \({\cal H}\) we choose \(-f\) and use the same tricks 
as before.  

\paragraph{\bf C. The region ${\cal N}_\eta$} Selecting accordingly we have that for 
$\chi_0=1-\chi_--\chi_+$ we have that
\[
 \int_{{\cal N}_\eta} \chi_0 w \lap{{\cal H}}f\leq C\eta^\ell 
\]
for suitable $\ell$.

\paragraph{Gluing patches together} Finally we have that for appropriate choice of $\eta$ relaitve to $\epsilon$ we have that: 
\[
  4\delta \olwv{\pm,\eta} wf -\eta -2\olwv{}
(2Hw^2+h(\kl f,\kl f)  \geq 2n\olwv{}wf 
\] 
 Therefore we have that for minimal hypersurfaces then $H=0$ and 
selecting \({\cal S}_1,{\cal S}_2\) so that 
$h(\kl f,\kl f)>0$ we conclude that:
\[
 \delta > \frac{n}{2} 
\]

 \end{document}